\documentclass[a4paper, 12pt]{article}

\usepackage{amsmath, amsfonts, amssymb}
\usepackage{verbatim}
\usepackage[british]{babel}
\usepackage{graphicx}
\usepackage{color}

\usepackage{amsthm}

\newtheorem {proposition}{Proposition}
\newtheorem {theorem}{Theorem}
\newtheorem{remark}{Remark}

\newcommand{\R}{{\mathbb{R}}}

\newcommand{\Z}{\mathbb{Z}}
\newcommand{\N}{\mathbb{N}}

\newcommand{\eps}{\varepsilon}

\newcommand{\G}{\mathsf{L}}
\renewcommand{\S}{\mathsf{Supp}}
\renewcommand{\L}{\Lambda}
\newcommand{\ex}{e^{\left(x\right)}}

\newcommand{\bs}[1]{\boldsymbol{#1}}%
\newcommand{\A}{\tilde{A}}

\date{}

\begin{document}

\title{A note on scaling  limits for  truncated birth-and-death processes with interaction }

\author{
Vadim Shcherbakov\footnote{Department of Mathematics,
 Royal Holloway,  University of London.
 Email address: vadim.shcherbakov@rhul.ac.uk
}
\, and 
Anatoly Yambartsev\footnote{Department of Statistics,
 University of Sao Paulo.
 Email address: yambar@gmail.com
}
}


\maketitle


\begin{abstract}
{\small 
In this note we consider  a Markov chain formed by  a finite system of   interacting birth-and-death processes on a finite state space.
 We study an  asymptotic behaviour of the Markov chain  as its state space becomes large. 
In particular, we show that the  appropriately  scaled Markov chain  converges  to a diffusion process,
and derive conditions for  existence of the stationary distribution of  the limit diffusion process
in special  cases.}
\end{abstract}

\section{The model}
\label{model}

Many real life systems are multicomponent, where 
the evolution  of an isolated  single component  is relatively simple, but the  presence of an interaction  
affects both the individual behaviour of a component and the collective behaviour.  
Also, the time evolution of many real life systems   
can be often described in terms of certain birth and death events.
Models of interacting  birth-and-death processes on integers  
 provide a flexible  mathematical framework for modelling such systems
(e.g. see \cite{VS1}, \cite{VS2},  \cite{Thai},  \cite{Triolo} and references therein) that appear in biology, physics, queueing 
and other applications. Frequently, there are natural limitations on the system 
size (e.g. limited resources in biological systems, restrictions of the queue length in queueing etc.). On the other hand, a state space
of finite system  can be large. 
This motivation stimulates interest to finding an adequate asymptotic description of finite but large real life systems with the features described above.   

This note concerns scaling limits for a stochastic model that describes 
a finite  continuous Markov chain formed by interacting birth-and-death processes confined to a finite set.
The system components (spins) are  labelled by vertices of a finite connected graph and evolve  
subject to a local  interaction determined by the graph. The interaction  not necessarily symmetric.
It should be noted  that spin  models with  asymmetric interactions have  recently been introduced 
 for modelling interaction in biological systems  (see~\cite{Fontes} for  detailed explanations and references therein).
We are interested in the asymptotic behaviour of the  Markov chain  as the  range of  possible values  of  its 
components   becomes large. We show that the appropriately scaled Markov chain 
can  be approximated  either  a by  diffusion process, or by a  deterministic process depending on the scaling.  
Such scaling limits  are widely used  in queueing (e.g., see \cite{Kurtz},  \cite{Glynn},   \cite{Stone}), 
 interacting particle  systems  (e.g., see \cite{Velenik}, \cite{Landim} and references therein) and in many other applications (e.g. in finance,   \cite{Lykov} and  \cite{Muzychka}).

  Let us describe the model.   Let $\Lambda$ be   a finite connected  graph.
Given   integers $l\geq 0$ and $r>0$  define  $\Omega_{\Lambda, l, r}=\{-l,\ldots, r\}^{\L}$.
Denote by  $\xi_x,\, x\in \L$, components of $\xi\in \Omega_{\Lambda, l,r}$ and call them spins.
We write  $x\sim y$ to denote that vertices $x, y\in \L$ are  adjacent,  
and  $x\nsim y$ if they are not. Vertices  $x$ and $y$  are   called neighbours, if $x\sim y$.
By convention, $x\sim x$ for all $x\in \L$.
A matrix  $A=(\alpha_{xy})_{x,y\in \L}$ is called an interaction matrix, if  $\alpha_{xy}=0$, whenever 
 $x\nsim y$. It is easy to see  that 
 any linear combination of interaction matrices is an interaction matrix.
Given  two interaction matrices $A_b$ and $A_d$  consider 
a continuous time birth-and-death Markov chain $\xi(t)\in \Omega_{\Lambda, l,r}$ evolving as follows.
 Given that  $\xi(t)=\xi$  a  spin $\xi_x<r$ increases by $1$ at the rate  
$e^{b(x, \xi)}$, and a spin  $\xi_x>-l$ decreases by $1$ at the rate 
$e^{d(x, \xi)}$, where  $b(x, \xi)=(A_b\xi)_x$ and $d(x, \xi)=(A_d\xi)_x$.

It is easy to see that various types of interaction between  spins can be modelled   by choosing appropriate 
interaction matrices $A_b$ and $A_d$. For example, if these matrices 
are diagonal then components of the Markov chain are independent truncated birth-and-death processes.
In general,  matrices $A_b$ and $A_d$   are not 
symmetric.

A  variant of this Markov chain was considered in~\cite{VS1},  where 
 transition rates  were specified 
by   interaction matrices  $A_b=A=(\alpha_{xy})_{x,y\in \L}$,  such that $\alpha_{xy}\equiv const$, 
and  $A_d\equiv 0$,  and   graph $\L$  was  a   $d$-dimensional lattice cube.
It was shown that a stationary distribution of the Markov chain 
converges to a Gibbs measure,  as $\L$ expands to the whole lattice,   
  and an  occupied site percolation problem  was solved for the  limit distribution. 
 The long term behaviour of the   Markov chain with non-negative and 
 unbounded components (formally obtained by setting  $l=0$ and $r=\infty$)  
was studied  in~\cite{VS2}. The transition rates in ~\cite{VS2} 
 were specified  by interaction matrices  $A_d\equiv 0$ (as in~\cite{VS1}) and $A_b=\alpha E +\beta I_{\L}$,
where  $\alpha, \beta\in \R$, $I_{\L}$ is  the  incidence matrix of graph $\L$ and $E$ is the unit  matrix.  
 The main goal  in~\cite{VS2}
 was to determine how the  long term behaviour of the Markov chain depends on both the transition parameters 
and the structure of  the underlying graph.  
The model in this note is  somewhat intermediate  between the models   in~\cite{VS1} and  in~\cite{VS2}.
Namely, the underlying graph is fixed and we study  the asymptotic behaviour of the Markov chain, as  the finite 
range of the spin values expands. 
If we formally equate $l=-\infty$ and $r=\infty$, then the corresponding countable continuous Markov chain can be 
explosive (this depends on both  the interaction matrices and graph $\L$).
On the other hand, if  we stretch the range of  the  spin   values 
and   simultaneously change  the model parameters, then the Markov chain can converge
to  a non-trivial limit  under an appropriate space-time scaling.

  The diffusion limit   (Theorem \ref{T2}) is of a particular interest as it can be interpreted in terms of 
a system of interacting one-dimensional diffusions.   
In some  cases   these   one-dimensional diffusions are given by  famous Ornstein-Uhlenbeck 
processes. It is known that a single Ornstein-Uhlenbeck process on the line is positive recurrent with 
a stationary distribution given by a Gaussian probability density. 
Presence of interaction can significantly change the collective behaviour 
of the system. Namely, the  system can become transient.
This effect  depends on both interacting matrices and the structure of  graph $\L$
as demonstrated by examples in Section~\ref{rev}.
  
We also formulate  an analogue of the law of large numbers (Theorem~\ref{T1}), where  the limit process is  
described by system of ordinary differential equations.

\section{Scaling limits}

%

Given  interaction matrices $A_b$ and $A_d$ 
consider  a Markov  process $\bs{u}(t)=\{u_x(t),\, x\in \L\}\in \R^{\L}$ 
which   is  a solution of the following system of stochastic differential equations
\begin{align}
\label{Lim}
du_x(t)&=\left(b(x, \bs{u}(t))-d(x, \bs{u}(t))\right)dt + \sqrt{2}dW_x(t),\quad x\in \L,\\
u_x(0)&=u_x,\, x\in \L,\nonumber
\end{align}
where
 $b(x, \xi)=(A_b\xi)_x$ and $d(x, \xi)=(A_d\xi)_x$,
 and  $W_x(t)\in \R,\, x\in \L, $  are independent one-dimensional  standard Brownian motions.
Define the following matrix 
 \begin{equation}
\label{A}
A=A_b-A_d=(\alpha_{xy})_{x,y\in \L}.
\end{equation}
Equations (\ref{Lim}) can be now  rewritten as follows
\begin{align*}
du_x(t)&=\left(\alpha_{xx}u_x(t)+\sum_{y\neq x}\alpha_{xy}u_y(t)\right)dt + \sqrt{2}dW_x(t),\quad x\in \L,\\
u_x(0)&=u_x,\, x\in \L,\nonumber
\end{align*}
or,  in the following  vector form 
\begin{align}
\label{Linear1}
d\bs{u}(t)&=A\bs{u}(t)dt +\sqrt{2}dW(t),\\
\bs{u}(0)&=\bs{u}\in \R^{\L}, \nonumber
\end{align}
where $W(t)=\{W_x(t)\in \R,\, x\in \L\}\in \R^{\L} $.

\begin{remark}
\label{OUrem}
{\rm 
Note that if diagonal elements of matrix  $A$ are  negative, i.e.
$\alpha_{xx}<0$ for all $x\in \L$, then
diffusion  process $\bs{u}(t)$ can be interpreted as   a   system of locally interacting Ornstein-Uhlenbeck processes
with individual drifts $\alpha_{xx}$, $x\in \L$ and constant diffusion coefficients equal to $\sqrt{2}$.
}
\end{remark}

\begin{theorem}
\label{T1} 
Given  interaction matrices $A_b$ and $A_d$, and sequences of positive numbers  
$\eps_n, l_n, r_n,\, n\in \N,$ consider a sequence of 
 Markov chains  $\xi^{\left(n\right)}(t)\in \Omega_{\L, l_n, r_n},\, n\in \N,$
whose transition  rates are specified by interaction matrices  $\eps^2_nA_b$ and $\eps^2_nA_d$, 
Suppose that  
\begin{equation}
\label{ln-rn1}
\lim_{n\to \infty}\eps_n=0, \,  \lim_{n\to \infty}l_n\eps_n= \lim_{n\to \infty}r_n\eps_n=\infty,\, 
\lim\limits_{n\to \infty}\eps_n\xi^{\left(n\right)}\left(0\right)=u\in\R^{\L}.
\end{equation}
Under these conditions the sequence 
of rescaled Markov chains
$\eps_n\xi^{\left(n\right)}\left(t\eps_n^{-2}\right)$
converges as $n\to \infty$ to a Markov process 
 $\bs{u}(t)\in \R^{\Lambda}$  which is a unique solution of equation~(\ref{Lim}) with initial condition $\bs{u}(0)=u$.
The convergence is understood in a sense of the weak convergence of the corresponding semigroups.
\end{theorem}
The proof of Theorem~\ref{T1} is given in Section~\ref{proofT1}.

\begin{remark}
\label{R1}
{\rm 
 Both the  Markov chain and the diffusion limit  are Feller Markov processes.
For Feller Markov processes weak convergence is equivalent to convergence of the corresponding semigroups
(Theorem 2.5, Chapter 4, \cite{Kurtz}). }
\end{remark}
The following theorem is an analogue of the law of large numbers.
\begin{theorem}
\label{T2}
Given interaction matrices $A_b$ and $A_d$, and sequences of positive numbers  
$\eps_n, l_n, r_n,\, n\in \N,$ consider a sequence of 
 Markov chains  
$\xi^{\left(n\right)}(t)\in \Omega_{\L, l_n, r_n},\, n\in \N,$ 
whose transition rates are specified by interaction matrices 
$\eps_nA_b$ and $\eps_nA_d$.
Suppose that 
\begin{equation}
\label{ln-rn2}
\lim_{n\to \infty}\eps_n=0, \,  \lim_{n\to \infty}l_n\eps_n= \lim_{n\to \infty}r_n\eps_n=\infty,
\, \lim\limits_{n\to \infty}\eps_n\xi^{\left(n\right)}\left(0\right)=u \in\R^{\L}.
\end{equation}
Under these conditions  for every $t\geq 0$
\begin{equation}\label{conv1}
\lim_{n\to\infty} \sup_{s\le t} \bigl| \eps_n\xi^{\left(n\right)}\left(t\eps_n^{-1}\right) - \gamma(t) \bigr| =0, \mbox{ a.s. }
\end{equation}
where deterministic process  $\gamma(t)=\bigl( \gamma_x(t), x\in\Lambda \bigr)$ solves the following system of 
non-linear differential equations
\begin{equation}\label{sys1}
\dot{\gamma}_x(t) = e^{b(x, \gamma(t)) } - e^{d(x, \gamma(t)) } ,\,\, x\in \L,
\end{equation}
with initial  conditions $\gamma(0)=u$.
\end{theorem}

\begin{remark}
\label{R2}
{\rm 
Note that  assumptions 
(\ref{ln-rn1}) and~(\ref{ln-rn2}) of Theorem~\ref{T1} and Theorem~\ref{T2} respectively
can be modified  as follows
$$\lim\limits_{n\to \infty}\eps_n=0,\,  \lim_{n\to \infty}l_n\eps_n=a,\,  \lim_{n\to \infty}r_n\eps_n=b,\,
\lim\limits_{n\to \infty}\eps_n\xi^{\left(n\right)}\left(0\right)=\bs{u}(0)\in [-a, b]^{\L},$$
where both  $a$ and $b$ can be either  finite, or infinite. 
}
\end{remark}

\section{ Invariant measures and the diffusion limit  in the reversible  case}
\label{rev}

In this section we consider invariant measures of both the Markov chain and the diffusion limit in  the reversible case.
 Let us assume, throughout the section, that  matrix $A=A_b-A_d$$=(\alpha_{xy})_{x,y\in \L}$ is  symmetric  (the symmetric case).
In the symmetric case the Markov   chain  is reversible with  the following stationary distribution 
\begin{equation}
\label{measure}
\mu_{\L}(\xi)=Z_{\L}^{-1}e^{
\frac{1}{2}\sum\limits_{x}\alpha_{xx}\xi_x(\xi_x-1)+ \sum\limits_{x\sim y}\alpha_{xy}\xi_x\xi_y},\,\quad \xi\in \Omega_{\L, l, r},
\end{equation}
where  
$$
Z_{\L}=\sum\limits_{\xi\in \Omega_{\L, l, r}}e^{
\frac{1}{2}\sum\limits_{x}\alpha_{xx}\xi_x(\xi_x-1)+ \sum\limits_{ x\sim y}\alpha_{xy}\xi_x\xi_y}.
$$
Indeed, it is easy to see that  measure (\ref{measure}) satisfies the following detailed balance equation 
 $$e^{(A_b\xi)_x}\mu_{\L}\left(\xi\right)
=\mu_{\L}\left(\xi+e^{\left(x\right)}\right)e^{(A_d\xi)_x},$$
or, equivalently, 
 $$e^{((A_b-A_d)\xi)_x}\mu_{\L}\left(\xi\right)
=\mu_{\L}\left(\xi+e^{\left(x\right)}\right),$$
where addition of configuration 
is understood component-wise and  $e^{\left(x\right)}\in \Omega_{\L, l, r}$ 
is the configuration such that $e^{\left(x\right)}_y=0$ if $y\neq x$, and  $e^{\left(x\right)}_x=1$.
 In vector notation
$
\mu_{\L}(\xi)=Z_{\L}^{-1}e^{
\frac{1}{2}\left[(A\xi, \xi)-(\alpha, \xi)\right]} $
and
$
Z_{\L}=\sum_{\xi\in \Omega_{\L, l, r}}e^{\frac{1}{2}\left[(A\xi, \xi)-(\alpha, \xi)\right]},
$
where  $\xi\in \Omega_{\L, l, r}$, $\alpha=(\alpha_{xx},\, x\in \L)$ is a vector formed by diagonal elements of matrix $A$, and
$(\xi', \xi'')$ is  the Euclidean scalar product of vectors $\xi', \xi''\in \Omega_{\L, l, r}$ (considered as  elements of $\R^{\L}$).

\begin{remark}
{\rm 
If $l=0, r=1$ and $\alpha_{xy}\equiv  const$, then probability distribution~(\ref{measure}) corresponds to a particular case of the  celebrated   Ising model}. 
\end{remark}

It is easy to see that under conditions of Theorem~\ref{T1} if  a 
 sequence of states $\xi^{\left(n\right)}\in \Omega_{\L, l_n, r_n},\, n\in \N,$ is such that 
$\eps_n\xi^{\left(n\right)}\to u\in \R^{\L}$ as $n\to \infty$, 
then 
\begin{equation}
\label{conv}
e^{\frac{1}{2}\left[(\eps_n^2A\xi^{\left(n\right)}, \xi^{\left(n\right)})-
(\eps_n^2\alpha, \xi^{\left(n\right)})\right]}\to e^{\frac{1}{2}(A u, u)}, \quad u\in \R^{\L},
\end{equation}
where the function in the right side
should be a density of an invariant measure of the diffusion limit. If this density is integrable, then the properly normalised 
invariant measure is the stationary 
distribution of the diffusion limit and, hence, the latter is positive recurrent.
In what follows we are going to consider conditions for existence of the stationary distribution of the diffusion limit and, therefore, existence 
of the stationary distribution. 
  
Note  first that   a unique strong solution 
of equation~(\ref{Linear1}), regardless of symmetry of $A$,   is given by the following formula  (e.g. Section 5.6 in~\cite{Karatzas})
\begin{equation}
\label{Solution}
\bs{u}(t)=e^{At}\bs{u}(0)+\sqrt{2}\int\limits_{0}^{t}e^{A(t-s)}dW_s,
\end{equation}
where $e^{At}=\sum_{n=0}^{\infty}\frac{t^n}{n!}A^n$.
 If  all eigenvalues of $A$ have negative real parts,
then  the process is positive recurrent and its stationary distribution 
is a zero mean Gaussian distribution with  the covariance  that can be expressed in terms of matrix $A$ 
(e.g. Theorem 6.7,  Section 5.6, \cite{Karatzas}). 

If $A$ is symmetric, then  equation (\ref{Linear1}) can be written in the following gradient form 
\begin{equation}
\label{Grad}
d\bs{u}(t)=\frac{1}{2}\nabla (Au(t),u(t))dt +\sqrt{2}dW(t),
\end{equation}
which is a particular case of  the Langevin equation. It follows from the general theory of Langevin equations that 
the diffusion process $\bs{u}(t)$  is reversible and the following 
function  (the same as in~(\ref{conv}), as it should be)
$$
e^{\frac{1}{2}(Au, u)}=e^{-\frac{1}{2}(\tilde A u, u)}=
e^{\frac{1}{2}\sum_{x}\alpha_{xx}u_x^2+\sum_{x\sim y}\alpha_{xy}u_xu_y}
,\,\, u=\{u_x,  x\in \Lambda\}\in \R^{\L},
$$
 is a density of an invariant  measure of the process.
This density is integrable if and only if matrix  $\tilde{A}=-A$ is positive definite,
in which  case   a stationary distribution of diffusion process $\bs{u}(t)$ is
a multivariate normal distribution with zero mean and covariance matrix $(\A)^{-1}$. 
In the rest of the section we are going to obtain conditions of positive  definiteness 
of symmetric $-A$ in special cases.

We start with noticing that for matrix $-A$ to be positive the  diagonal elements must be positive, which means that 
bounds 
$\alpha_{xx}<0$ must hold for all $x\in \L$. Note that in this case
the diffusion limit can be interpreted in terms of interacting Ornstein-Uhlenbeck processes (see Remark~\ref{OUrem}).

It is known from algebra, that  a  diagonally dominant symmetric matrix with positive elements on the main diagonal is positive definite. 
This  implies in our case  that if $\alpha_{xx}<0$ and $\alpha_{xx}+\sum_{y\sim x}|\alpha_{xy}|<0$ for all $x\in \L$, then
matrix $\tilde{A}$ is positive definite. In turn, this fact  implies  the following proposition.
\begin{proposition}
\label{P1}
Let   $\L$ be   an arbitrary finite  connected graph  and let $I_{\L}$ be  the  incidence matrix of $\L$.
If $A=\alpha E + \beta I_{\L}$, 
where $\alpha, \beta\in \R$,  and $E$  is   the unit matrix, 
and 
$\alpha<0,\, \alpha+|\beta|\max_{x\in \L}\nu(x)<0$, 
where $\nu(x)$  is the degree of vertex $x$ (i.e.  the number of  edges incident to the vertex),
then $\A$ is positive definite.
\end{proposition}
Inequality $|\beta|\max_{x\in \L}\nu(x)<-\alpha$  in the above proposition 
means that  interaction (specified by parameter $\beta$ and graph $\L$) is sufficiently small, so that 
the collective behaviour of the system of interacting Ornstein-Uhlenbeck processes (with individual drifts equal $\alpha$) 
is still positive recurrent.
 
An additional information about graph $\L$  allows  to improve this result, which we are 
going to demonstrate in the case of the following  graphs.
\begin{enumerate}
 \item A constant vertex degree graph is a graph such that  
 $\nu(x)\equiv \nu$, for some integer $\nu>0$, where $\nu(x)$ is the degree of vertex $x$.
\item  A star  graph with $m+1$ vertices  is a graph with central vertex  $x$ and its neighbouring vertices   
$y_1, \ldots, y_m$, i.e.  $x\sim y_i,\, i=1,\ldots,m,$ so that  $x$ is the only neighbour  for  each of $y_i,\, i=1,\ldots,m$.  
\item A  unary tree of length $n+2$, where $n\in \Z_{+}$, 
 is a graph which vertices can be enumerated by natural numbers $1,\ldots, n+2$, and such that 
$1\sim 2 \sim \cdots \sim n+1\sim n+2$. 
\end{enumerate}
  The following theorem gives criteria for positive definiteness of matrix $\A$ and, hence, for positive recurrence 
of the corresponding diffusion limit  in the case of these graphs.
\begin{theorem}
\label{P2}
Suppose that $A=\alpha E + \beta I_{\L}$, 
where $\alpha, \beta\in \R$, $E$ is the unit matrix and $I_{\L}$ is the incidence matrix of $\L$.

1) 
If $\L$ is a graph  with constant vertex degree $\nu(x)\equiv \nu\geq 1$, then 
 $\A$ is positive definite if and only if 
$\alpha<0,\, \alpha+|\beta|\nu<0$.

2) If   $\L$ is a star-like  graph with $(m+1), m\geq 2$, vertices,  then 
$\A$ is positive definite 
 if and only if $\alpha<0,\, \alpha+|\beta|\sqrt{m}<0$.

3) If $\L$ is a unary tree of length $n+2$, where $n\in \Z_{+}$, then 
 then $\A$ is positive definite 
 if and only if $\alpha<0, \alpha+2\beta\cos\left(\frac{\pi}{n+3}\right)<0$.
\end{theorem}

\paragraph{{\it Proof of  Part 1) of Theorem~\ref{P2}.}}  The "if'' part of the statement is implied by Proposition~\ref{P1}.
To show that the sufficient condition is also a necessary one it suffices to  notice 
that all eigenvalues of matrix $\tilde A$  lie, by the Gershgorin circle theorem,
within the closed interval $[-\alpha-|\beta|\nu, -\alpha+|\beta|\nu]$.

\paragraph{{\it Proof of Part 2) of Theorem~\ref{P2}.} } Let  $\L$ be  a star  graph 
 with a  central vertex  $x$ and its neighbouring vertexes  $y_1, \ldots, y_m$, i.e.  $x\sim y_i,\, i=1,\ldots,m,$ and $x$ is the only neighbour  for  each of $y_i,\, i=1,\ldots,m$.  
  Denote by $D_m(\mu)$ the characteristic polynomial 
of matrix $\A$ corresponding to the graph. 
It was shown in~\cite{VS2} that 
$$D_m(\mu)=(-\alpha-\mu)^{m-1}(-\alpha-\beta\sqrt{m}-\mu)(-\alpha+\beta\sqrt{m}-\mu),$$
so that   $-\alpha>0$ is the matrix eigenvalue of order $m-1$ and  $-\alpha\pm \beta\sqrt{m}>0$
are two remaining  eigenvalues, each of order $1$ and the result follows.
 
\paragraph{{\it Proof of Part 3) of Theorem~\ref{P2}.}}
If $n=0$, then this is the simplest case of  a constant degree graph (see Part 1)).
If $n=1$, then this is the simplest case of  a star graph (see Part 2)). 
In what follows we  assume that $n\geq 2$.
It is easy to  see that matrix $\A$ is the following  tridiagonal symmetric Toeplitz matrix 
$$
\A=\begin{bmatrix}
-\alpha & -\beta &  & & & & 0  \\
-\beta & -\alpha & -\beta &  & &  \\
 & -\beta & \cdot & \cdot & & & \\
& & \cdot & \cdot & \cdot &  & \\
& & & \cdot & \cdot & \cdot &  \\
& & & & \cdot & \cdot & -\beta\\
0 & & & & & -\beta & -\alpha
\end{bmatrix}_{(n+2)\times (n+2)}
$$
The well known  results for  tridiagonal symmetric Toeplitz matrices 
yield  that  eigenvalues of matrix $\tilde A$ are simple and  given by the following equations
$$\lambda_k=-\alpha-2\beta\cos\left(\frac{k\pi}{n+3}\right),\, k=1,\ldots, n+2,$$
where  $\lambda_1=-\alpha-2\beta\cos\left(\frac{\pi}{n+3}\right)$ is the minimal eigenvalue. This finishes the proof of Part 3)
of  the theorem.

\section{Proof of Theorem \ref{T1}}
\label{proofT1}
The proof consists in showing that  a sequence of  generators of rescaled Markov chains converges in a certain sense (explained below)
to the generator of the limit diffusion process   \cite{Kurtz}.
 Solution of~(\ref{Lim})  is a  Feller Markov process and its 
generator is 
\begin{equation}
\label{G}
\G f(u)=\sum\limits_{x\in \Lambda}f''_{xx}(u) + 
\sum\limits_{x\in \Lambda}\left(b(x, u)-d(x, u)\right)f'_x(u), \quad  
u\in \R^{\L},
\end{equation}
where we denoted 
$f''_{xx}(u)=\frac{\partial^2 f(u)}{\partial u^2_x}$ and $f_x'(u)=\frac{\partial f(u)}{\partial u_x}$
for the ease of notation.
Generator $\G$  is defined on  twice continuously  differentiable functions that  vanish at infinity together 
with their first and second order derivatives.
By Theorem 6.1, Chapter 1, \cite{Kurtz},   in order to prove the convergence of 
semigroups it is sufficient to prove 
that for each  $f\in D(\G)$, where $D(\G)$ is a  core  of the limit generator
$\G$, there exists 
$f_{n}\in B(\eps_n\Omega_{\L, n}),\, n\geq 1,$ such that   
$$\sup_{\xi^{\left(n\right)}\in \Omega_{\L, n}}|f_n(\eps_n\xi^{\left(n\right)})-f(\eps_n\xi^{\left(n\right)})|\to 0$$
and 
$$\sup_{\xi^{\left(n\right)}\in \Omega_{\L, n}}|\G_nf_n(\eps_n\xi^{\left(n\right)})-\G f(\eps_n\xi^{\left(n\right)})|\to 0$$
as $n\to \infty$,
 where $ \G_n$ is the generator of  the rescaled Markov chain $\eps_n\xi^{\left(n\right)}(\eps_n^{-2}t)$.
 Theorem 2.5,  Chapter 8, \cite{Kurtz}, yields 
 that the set $C_c^{\infty}(\R^d)$ of infinitely differentiable 
functions with a compact support is   a core for  generator~(\ref{G}).

Denote  $ \Omega_{\L, n}= \Omega_{\L, l_n, r_n}$ to ease notation.
Let  $\xi^{\left(n\right)}(t)\in \Omega_{\L, n},\, n\in \Z_{+},$ 
be the Markov chain 
whose transition  rates are specified by interaction matrices  $\eps^2_nA_b$ and $\eps^2_nA_d$.
 Given  $f:C_c^{\infty}(\R^{\L})\to \R$ 
and $n$ define $f_n: \Omega_{\L, n}\to \R$ as $f_n=f|_{\Omega_{\L, n}}$, i.e. as a restriction of $f$ on $\Omega_{\L, n}$.
If $\G_n$ is the   generator of the rescaled Markov chain $\eps_n\xi^{\left(n\right)}(\eps_n^{-2}t)$, then 
we have   that 
\begin{align*}
\G_{n}f_n(\eps_n\xi)&=\eps_n^{-2}\sum\limits_{x\in \Lambda}\left(f(\eps_n(\xi+\ex))-
f(\eps_n\xi)\right)e^{b_n(x, \xi)}{\bf 1}_{\{\xi_x<r_n\}}\\
&+\eps_n^{-2}\sum\limits_{x\in \Lambda}\left(f\left(\eps_n(\xi-\ex)\right)-f(\eps_n\xi)\right)e^{d_n(x, \xi)}{\bf 1}_{\{\xi_x>-l_n\}},
\end{align*}
where we denoted $b_n(x, \xi)=\eps_n^2b(x, \xi)(=\eps_n^2(A_b\xi)_x)$, $d_n(x, \xi)=\eps_n^2d(x, \xi)(=\eps^2_n(A_d\xi)_x)$ and 
 ${\bf 1}_{B}$ is an indicator of  set  $B$.

Suppose that $n$ is   sufficiently large so that 
support  $\S(f)$ of $f$ is covered by $[-\eps_n l_n, \eps_n r_n]^{\L}$. Then 
\begin{align*}
\G_{n}f_n(\eps_n\xi)&=\eps_n^{-2}\sum\limits_{x\in \Lambda}\left(f(\eps_n(\xi+\ex))-
f(\eps_n\xi)\right)e^{b_n(x, \xi)}\\
&+\eps_n^{-2}\sum\limits_{x\in \Lambda}\left(f\left(\eps_n(\xi-\ex)\right)-f(\eps_n\xi)\right)e^{d_n(x, \xi)},
\end{align*}
in other words, we can remove the indicators ${\bf 1}_{\{\xi_x<r_n\}}$ and  ${\bf 1}_{\{\xi_x>-l_n\}}$.
 Note that due to linearity we have that  $b_n(x, \xi)=\eps_nb(x, \eps_n\xi)$ and  $d_n(x, \xi)=\eps_nd(x, \eps_n\xi)$. 
If  a sequence of states $\xi^{\left(n\right)}\in \Omega_{\L, n},\, n\geq 1,$ is such that 
$\eps_n\xi^{\left(n\right)}_{x}\to u_x$ for every $x\in \L$  as $n\to \infty$, then by Taylor's formula 
with the reminder term we can write that
\begin{align}
\label{tay1}
e^{\eps_nb(x, \eps_n\xi)}&=(1+\eps_nb(x,  u))+R_{n,1}(x, u),\\
\label{tay11}
e^{\eps_nd(x,\eps_n\xi)}&=(1+\eps_nd(x,  u))+R_{n,2}(x, u)
\end{align}
where  $|R_{n,i}(x,u)|<C_i\eps_n^2$, for some $C_i=C_i(f)$, $i=1,2$, 
as  $n\to \infty$. Also,   
\begin{equation}
\label{plus}
f(\eps_n(\xi+\ex)) - f(\eps_n\xi)=f_x'(u)\eps_n+
\frac{1}{2} f''_{xx}(u)\eps_n^2+R_{n, 3}(x, u),
\end{equation}
where $|R_{n, 3}(x, u)|\leq C_3\eps_n^3$, with some $C_3=C_3(f)$, and 
\begin{equation}
\label{min}
f\left(\eps_n(\xi-\ex)\right)-f(\eps_n\xi)=-f_x'(u)\eps_n+
f_{xx}''(u)\eps_n^2+R_{n, 4}(x, u),
\end{equation}
where $|R_{n, 4}(x, u)|\leq C_4\eps_n^3$, with some $C_4=C_4(f)$.
Equations~(\ref{tay1}), (\ref{tay11}),  
~(\ref{plus}) and~(\ref{min})  yield  that 
$$\G_{n}f_n(\eps_n\xi)=f_{xx}''(u)+(b(x, u)-d(x, u))f_x'(u)+J_{n}(x,u),
$$
where $|J_{n}(x,u)\leq C\eps_n^3$, $C=C(f)$.
Thus we have that 
$$
\G_{n}f(\eps_n\xi^{\left(n\right)})\to \sum\limits_{x\in \Lambda}\left(f_{xx}''(u)
+(b(x, u)+d(x, u))f_x'(u)\right),
$$
uniformly over $u\in \S(f)$.  Now, Theorem 6.1, Chapter 1, \cite{Kurtz} applies  and the convergence of semigroups 
follows.



\begin{thebibliography}{1}
\bibitem{Kurtz} Either. S.  and  Kurtz, T. (1985). {\it Markov processes: characterisation and convergence}.
John Wiley$\&$Sons.
\bibitem{Fontes}  Fernandez, R., Fontes, Luiz R., and Neves, J. (2009).
Density-profile  processes describing biological signaling  networks: almost sure 
convergence to deterministic trajectories.
{\it Journal of Statistical Physics}, {\bf 136}, pp.~875--901. 
\bibitem{Glynn} Glynn, P. and Ward, A. (2003).  
A diffusion approximation for a Markovian queue with reneging. {\it Queueing Systems},
{\bf 43}, pp.~103–128.
\bibitem{Velenik} Ioffe, D., Shlosman, S., and Velenik, Y. (2015). An invariance principle to Ferrari-Spohn diffusions.
{\it Communications in Mathematical Physics},   {\bf 336}, pp.~ 905--932.
\bibitem{Karatzas} Karatzas, I., and Shreve, S. (1991). {\it Brownian Motion and Stochastic Calculus}. Springer-Verlag, New York, 2nd Edition. 
\bibitem{Landim}
Kipnis, C. and Landim, C. (1999). {\it Scaling limits for interacting particle systems}. Springer-Verlag.
 \bibitem{Lykov} Lykov, A., Muzychka, S. and Vaninsky, K. (2012). Investor's sentiment in multi-agent model 
of the continuous auction. {\it arXiv:1208.3083v4.}
\bibitem{Muzychka} Muzychka, S., and Vaninsky, K. (2011). 
A class of nonlinear random walks related to the Ornstein-Uhlenbeck process.
{\it Markov Processes and Related Fields}, {\bf 17}, pp.~277--304.
\bibitem{VS1} Shcherbakov, V., and Yambartsev, A. (2012). On equilibrium distribution of a
reversible growth model. {\it Journal of Statistical Physics}, {\bf 148}, N1, pp.~53-66.
\bibitem{VS2} Shcherbakov, V., and Volkov, S. (2015). 
Long term  behaviour  of  locally interacting  birth-and-death processes. 
{\it Journal of Statistical Physics}, {\bf 158}, N1, pp.~132--157. 
\bibitem{Stone} Stone, C.  (1963). Limit theorems for random walks, birth and death processes, and diffusion
processes. {\it Illinois J. Math.},  {\bf 7}, pp.~638--660.
\bibitem{Thai}  Thia, M.-N. (2015).  Birth and Death process in mean field type interaction.
	{\it arXiv:1510.03238}.
\bibitem{Triolo} Triolo, L. (2005). Space Structures and Diﬀerent Scales for Many-Component Biosystems.
{\it Markov Processes and Related Fields}, {\bf 11}, pp.~389--404.


 \end{thebibliography}
\end{document}